\documentclass[12pt]{article}

\usepackage{amscd,amsmath, amssymb, fancyhdr, epsfig}
\usepackage[matrix,arrow]{xy}

\numberwithin{equation}{section}

\newcommand{\version}{version 1.0,\ \   Dec. 14, 2012}

\def\eqref#1{(\ref{#1})}

\newcommand{\Z}{{\Bbb Z}}
\newcommand{\C}{{\Bbb C}}

\newcommand{\R}{{\Bbb R}}

\renewcommand{\H}{{\Bbb H}}

\def\1{\sqrt{-1}\:}
\newcommand{\restrict}[1]{{\left|_{{\phantom{|}\!\!}_{#1}}\right.}}
\newcommand{\cntrct}                
{\hspace{2pt}\raisebox{1pt}{\text{$\lrcorner$}}\hspace{2pt}}

\makeatletter
\def\x@arrow{\DOTSB\Relbar}
\def\xlongequalsignfill@{\arrowfill@\x@arrow\Relbar\x@arrow}
\newcommand{\xlongequal}[2][]{%
        \ext@arrow 0099\xlongequalsignfill@{#1}{#2}}
\def\xlongrightarrowfill@{\arrowfill@\relbar\relbar\longrightarrow}
\newcommand{\xlongrightarrow}[2][]{%
        \ext@arrow 0099\xlongrightarrowfill@{#1}{#2}}
\makeatother

\newcommand{\calo}{{\cal O}}

\renewcommand{\phi}{\varphi}
\renewcommand{\epsilon}{\varepsilon}
\renewcommand{\geq}{\geqslant}
\renewcommand{\leq}{\leqslant}

\newcommand{\Sp}{\operatorname{Sp}}
\newcommand{\End}{\operatorname{End}}

\newcommand{\Id}{\operatorname{Id}}

\newcommand{\Hol}{\operatorname{Hol}}

\newcounter{Mycounter}[section]
\newcounter{lemma}[section]
\renewcommand{\thelemma}{{Lemma \thesection.\arabic{lemma}}}
\newcommand{\lemma}{%
    \setcounter{lemma}{\value{Mycounter}}
    \refstepcounter{lemma}
    \stepcounter{Mycounter}
    {\noindent \bf \thelemma:\ }}

\newcounter{claim}[section]
\newcounter{sublemma}[section]
\newcounter{corollary}[section]
\newcounter{theorem}[section]
\renewcommand{\thetheorem}{{Theorem \thesection.\arabic{theorem}}}
\newcommand{\theorem}{%
    \setcounter{theorem}{\value{Mycounter}}
    \refstepcounter{theorem}
    \stepcounter{Mycounter}
    {\noindent \bf \thetheorem:\ }}

\newcounter{conjecture}[section]
\newcounter{proposition}[section]
\newcounter{definition}[section]
\renewcommand{\thedefinition}
      {{Definition~\thesection.\arabic{definition}}}
\newcommand{\definition}{%
    \setcounter{definition}{\value{Mycounter}}
    \refstepcounter{definition}
    \stepcounter{Mycounter}
    {\noindent \bf \thedefinition:\ }}

\newcounter{example}[section]
\renewcommand{\theexample}{{Example \thesection.\arabic{example}}}
\newcommand{\example}{%
    \setcounter{example}{\value{Mycounter}}
    \refstepcounter{example}
    \stepcounter{Mycounter}
    {\noindent \bf \theexample:\ }}

\newcounter{remark}[section]
\renewcommand{\theremark}{{Remark \thesection.\arabic{remark}}}
\newcommand{\remark}{%
    \setcounter{remark}{\value{Mycounter}}
    \refstepcounter{remark}
    \stepcounter{Mycounter}
    {\noindent \bf \theremark:\ }}

\newcounter{problem}[section]
\newcounter{question}[section]
\makeatletter

\@addtoreset{equation}{section} \@addtoreset{footnote}{section} \makeatother

\def\blacksquare{\hbox{\vrule width 5pt height 5pt depth 0pt}}
\def\endproof{\blacksquare}

\begin{document}
\begin{center}
{\LARGE\bf
Holomorphic Lagrangian fibrations on hypercomplex manifolds\\[4mm]
}

Andrey Soldatenkov\footnote{Andrey Soldatenkov is partially supported by
AG Laboratory NRU-HSE, RF government grant, ag. 11.G34.31.0023},
Misha Verbitsky\footnote{Misha Verbitsky is partially supported by RFBR grant 10-01-93113-NCNIL-a,
RFBR grant 09-01-00242-a, Science Foundation of
the SU-HSE award No. 10-09-0015 and AG Laboratory NRU-HSE, 
RF government grant, ag. 11.G34.31.0023.}

\end{center}

{\small \hspace{0.10\linewidth}
\begin{minipage}[t]{0.85\linewidth}
{\bf Abstract} \\
A hypercomplex manifold is a manifold
equipped with a triple of complex structures
satisfying the quaternionic relations.
A holomorphic Lagrangian variety on a
hypercomplex manifold with trivial canonical bundle
is a holomorphic subvariety which is calibrated
by a form associated with the holomorphic volume form;
this notion is a generalization of the usual
holomorphic Lagrangian subvarieties
known in hyperk\"ahler geometry. An HKT 
(hyperk\"ahler with torsion) metric on 
a hypercomplex manifold is a metric
determined by a local potential, in a similar
way to the K\"ahler metric. We prove that
a base of a holomorphic Lagrangian fibration
is always K\"ahler, if its total space is HKT.
This is used to construct new examples
of hypercomplex manifolds which do not 
admit an HKT structure.

\end{minipage}
}

\tableofcontents


\section{Introduction}


\subsection{HKT metrics and $SL(n,{\Bbb H})$-structures}

Let $I,J,K$ be complex structures on a manifold $M$
satisfying the quaternionic relation $I\circ J=-J \circ I=K$.
Then $(M,I,J,K)$ is called {\bf a hypercomplex manifold}.
Hypercomplex manifolds are quaternionic analogues of
complex manifolds, in the same way as hyperk\"ahler
manifolds are quaternionic analogues of K\"ahler
manifolds. However, the theory of hypercomplex
manifolds is much richer in examples and (unlike
that of hyperk\"ahler manifolds) still not well understood.

For more details on hypercomplex and hyperk\"ahler
structures, please see Subsection \ref{_qD_Subsection_}.

The term {\it hypercomplex} is due to C. P. Boyer,
\cite{_Boyer_}, who classified compact hypercomplex
manifolds in $\dim_{\Bbb H}=1$. These are: K3 surfaces,
compact 2-dimensional complex tori, and Hopf surfaces.
However, in dimension $> 1$ there is no classification
apparent. 

The notion of hypercomplex structure
was considered as early as in 1955 by \cite{_Obata_},
who proved existence and uniqueness of a torsion-free
connection preserving a hypercomplex structure.
However, the first non-trivial examples
were obtained only in 1988 by physicists
Ph. Spindel, A. Sevrin, W. Troost, A. Van Proeyen 
(\cite{_SSTvP_}), and the intensive study of
hypercomplex structures began with the paper
\cite{_Joyce_} by D. Joyce, who classified
hypercomplex structures on homogeneous spaces,
and constructed one on each compact Lie group
multiplied by a compact torus of appropriate
dimension.

The main geometric tool allowing one to work with the
hypercomplex manifolds is a notion of an HKT
(``hyperk\"ahler with torsion'') metric,
due to physicists Howe and Papadopoulos
(\cite{_Howe_Papado_}). For a definition and
a discussion of HKT structure, please see 
Subsection \ref{_qD_Subsection_}.
HKT structures were much used in physics 
since mid-1990-ies (\cite{_GP1_}, \cite{_GP2_},
\cite{_GPS_}), and then mathematicians started
using HKT structures to study the hypercomplex
manifolds.

In \cite{_Verbitsky:HKT_}, the second named author
developed Hodge theory for HKT manifolds, and it turned
out to be quite useful. To illustrate the usability
of HKT metrics, let us quote the following result,
obtained in \cite{_V:hypercomplex_Kahler_}.
Let $(M,I,J,K)$ be a hypercomplex manifold with
$(M,I)$ a complex manifold of K\"ahler type.
Then $(M,I)$ is in fact hyperk\"ahler (that it,
admits a hyperk\"ahler structure). The idea 
of the proof is that any K\"ahler metric
on $(M,I)$ gives an HKT metric on $(M,I,J,K)$,
and this metric can be used to apply Hodge theory.
The interplay between K\"ahler geometry and
HKT geometry is quite extensive, and can be used
to study the geometry of hypercomplex manifolds
extensively.

In the present
paper, we obtain a result in a similar vein, constructing
a K\"ahler metric on a base of holomorphic Lagrangian 
fibration on an HKT manifold.

Of course, the notion of ``holomorphic Lagrangian
fibration'' is itself highly non-trivial, because
an HKT manifold is not necessarily holomorphically symplectic.
It was developed in \cite{_Gra_Verb_}, using the theory
of calibrations and the holonomy of the Obata connection.

The Obata connection is an immensely powerful tool,
known since 1950-ies (\cite{_Obata_}), but even the
most primitive invariants of Obata connection, such as
its holonomy, are hard to compute, and still unknown
in most cases.  When $(M,I,J,K)$ admits a hyperk\"ahler
metric, Obata connection is equal to the Levi-Civita
connection of the hyperk\"ahler metric and its holonomy
is $\Sp(n)$. The converse is also true: whenever
holonomy of a torsion-free connection is $Sp(n)$,
it is a Levi-Civita connection of a hyperk\"ahler manifold.

In a general situation, it's hard to say anything definite
about the holonomy $\Hol(\nabla)$ of the Obata connection. 
It is clear that $\Hol(\nabla)\subset GL({\Bbb H}, n)$,
but the equality is rarely realized. In fact, before
the publication of the paper \cite{_Soldatenkov:SU(3)_} in
2011, there was no single example of a hypercomplex
manifold with $\Hol(\nabla)= GL({\Bbb H}, n)$.
In \cite{_Soldatenkov:SU(3)_}, the first named author
provided an example of such a manifold, by proving that
$\Hol(\nabla)= GL({\Bbb H}, 2)$ for the homogeneous
hypercomplex structure on $SU(3)$, constructed
by D. Joyce in \cite{_Joyce_}.

However, there are many examples of hypercomplex manifolds
with $\Hol(\nabla)$ $\subsetneq GL({\Bbb H}, n)$.
Two biggest known families of hypercomplex manifolds are
the Joyce's homogeneous manifolds from \cite{_Joyce_},
and hypercomplex nilmanifolds (\cite{_Barberis+Dotti:1996_}).
As shown in \cite{_BDV:nilmanifolds_}, for all nilmanifolds,
the holonomy of the Obata connection lies in 
\[ 
SL({\Bbb H}, n)=[GL({\Bbb H}, n),GL({\Bbb H}, n)]\subset GL({\Bbb H}, n).
\]
A hypercomplex manifold with holonomy in $SL(n, {\Bbb H})$
is called {\bf an $SL(n, {\Bbb H})$-manifold}. Such manifold
might have holonomy group which is strictly less than
$SL(n, {\Bbb H})$; in fact, no example of a manifold
with $\Hol(\nabla)= SL({\Bbb H}, n)$ is known so far.

For any  $SL(n, {\Bbb H})$-manifold,  the Obata connection on 
the canonical bundle $K(M,I)=\Lambda^{2n,0}(M,I)$ of $(M,I)$ 
preserves a non-trivial section. This implies that $K(M,I)$
is holomorphically trivial (see \cite{_Verbitsky:canoni_}).
In the presence of an HKT metric, the converse is also true:
any compact hypercomplex manifold with 
holomorphically trivial canonical bundle
$K(M,I)$ and an HKT metric satisfies 
$\Hol(\nabla)\subset SL({\Bbb H}, n)$.
This result is obtained in \cite{_Verbitsky:canoni_}
using the Hodge theory for HKT manifolds
developed in \cite{_Verbitsky:HKT_}.

The Hodge-theoretic constructions
of \cite{_Verbitsky:HKT_} work for $SL(n, {\Bbb
  H})$-manifolds with HKT-structure especially well.
For such manifolds, one obtains Hodge-type decomposition
on the holomorphic cohomology bundle $H^*(\calo_{(M,I)})$.

\subsection{Calibrations on manifolds and Lagrangian fibrations}

Let $M$ be a Riemannian manifold. {\bf A calibration} on $M$
is a closed $k$-form $\eta$, such that $\eta(x_1, ...,
x_k)\leq 1$ for each orthonormal $k$-tuple $x_1, ...,
x_k\in TM$. A $k$-dimensional oriented subspace $V\subset T_x M$
is called {\bf calibrated} if the Riemannian volume of $V$
is equal to $\eta\restrict V$, and a subvariety
$Z\subset M$ of dimension $k$ is called {\bf calibrated} if its
singularities are of Hausdorff codimension $\geq 1$, and
all smooth tangent planes $T_z Z\subset T_z M$ are
calibrated.

The theory of calibrations has a long and distinguished
history, starting from the work of Harvey and Lawson 
\cite{_Harvey_Lawson:Calibrated_}. In \cite{_Gra_Verb_},
several families of calibrations were constructed
on hyperk\"ahler manifolds, using the quaternionic
linear algebra. Let $(M,I,J,K)$ be a hyperk\"ahler
manifold, and $\omega_I, \omega_J, \omega_K$ 
its K\"ahler forms. These forms generate an 
interesting commutative subalgebra in $\Lambda^*(M)$
(see \cite{_Gra_Verb_}, and \cite{_Verbitsky:trianalyt_} 
for structure results about this algebra). In
\cite{_Gra_Verb_}, Grantcharov and Verbitsky
discovered several new calibrations which are
expressed as polynomials of $\omega_I, \omega_J,
\omega_K$. One of these calibrations, denoted as $\Psi$
in the sequel (see \eqref{_eqn_Psi_}), is called 
{\bf holomorphic Lagrangian calibration}.
On a hyperk\"ahler manifold, it calibrates
holomorphic Lagrangian subvarieties in $(M,I)$.

It is surprising (and quite wonderful)
that this form remains
closed, even when the metric is not hyperk\"ahler,
provided that the Obata connection of $(M,I,J,K)$ belongs
to $SL(n, {\Bbb H})$. The hyperk\"ahler condition
can be weakened drastically: it suffices to assume
that the metric $g$ is quaternionic Hermitian, and
the Obata-parallel section of $K(M,I)$
has constant length with respect to $g$.
In this case $\Psi$ is closed. If rescaled properly, 
this form is a calibration, calibrating a
new class of complex subvarieties of $(M,I)$ 
called {\bf holomorphic  Lagrangian} (see
\cite{_Gra_Verb_}
and \ref{_holo_lagra_Definition_} for more details). 

In \cite{_SV:trianalytic_} it was shown that 
holomorphic Lagrangian subvarieties
of $SL(n, {\Bbb H})$-manifolds exist only
in a countable number of complex structure
of the form $L=aI+bJ+cK$ (such complex structures
are called {\bf induced by a hypercomplex structure},
see Subsection \ref{_qD_Subsection_}). In the present
paper, we study holomorphic Lagrangian fibrations
on hypercomplex manifolds. Such fibrations are often 
present in examples (\cite{_Gra_Verb_}; see also
Subsection \ref{_exa_Lagra_Subsection_}). 

The holomorphic Lagrangian fibrations are of significant
interest for mathematicians and physicists; for an early
survey of the subject, please see the paper
\cite{_Sawon_}, and for applications, see e.g.
\cite{_KRS_}. The non-K\"ahler version of this
geometry defined above should be even more
useful, because the number of examples is
much greater.

In the present paper, we consider HKT metrics on 
$SL(n, {\Bbb H})$-manifolds admitting holomorphic
Lagrangian fibrations. We prove that a base of such
fibration is always K\"ahler, if the fibration is 
smooth (\ref{_HKT=>Kahler_Theorem_}).
This allows us to construct new examples
of hypercomplex manifolds admitting no HKT metrics.

Originally, it was conjectured that any
compact hypercomplex manifold admits
an HKT metric. Examples of hypercomplex
nilmanifolds not admitting HKT metrics
were constructed by Fino and Grantcharov
(\cite{_Fino_Gra_}). Since then, hypercomplex
nilmanifolds admitting HKT metrics were
classified completely in \cite{_BDV:nilmanifolds_}.
It was found that in fact most of hypercomplex
nilmanifolds are not HKT. However, no other
hypercomplex manifolds not admitting an
HKT metrics were known before now.

\hfill

{\bf Acknowledgements.} The work on this paper began during
the visit to the University of Science and Technology of China.
The authors would like to express their gratitude to USTC
for the hospitality an to Prof. Xiuxiong Chen for the invitation.


\section{Preliminaries}


\subsection[Hypercomplex manifolds and quaternionic 
Dolbeault complex]{Hypercomplex manifolds \\ and 
quaternionic Dolbeault complex}
\label{_qD_Subsection_}

A $C^{\infty}$-manifold $M$ is called {\bf hypercomplex}
if $M$ is equipped with complex structures $I$, $J$, $K$,
that satisfy the quaternionic relations $$IJ = -JI =K.$$

It was shown in \cite{_Obata_} that any hypercomplex
manifold admits a unique torsion-free connection that preserves
$I, J$ and $K$. This connection is called
{\bf the Obata connection}. Note that the holonomy of the
Obata connection is a subgroup of $GL(n, \H)$.

Any almost-complex structure which is preserved by a torsion-free
connection is integrable (this follows from Newlander-Nirenberg
theorem). Therefore, for any $a, b, c\in \R$ with $a^2 + b^2 + c^2 = 1$,
the almost-complex structure $L = aI + bJ + cK$ is integrable.
We denote by $(M,L)$ the corresponding complex manifold.
Such complex structures are called {\bf induced by
quaternions}.

Let $(M, I, J, K)$ be a hypercomplex manifold of real dimension $4n$.
The hypercomplex structure induces the action of $SU(2)$ on all tensor
bundles over $M$. Recall that any irreducible complex representation
of $SU(2)$ is of the form $S^k U$, where $U$ is the standard
two-dimensional representation and $S^k$ denotes the symmetric power.
We will refer to any representation
of the form $\left(S^k(U)\right)^{\oplus m}$ (for arbitrary $m$)
as a {\bf weight k} representation of $SU(2)$.

We make the following observation about the weight decomposition
of the exterior algebra $\Lambda_\C^* M = \Lambda^* M \otimes_\R \C$.
First, note that $\Lambda^1_\C M$ is an $SU(2)$ representation
of weight one. It follows from the Clebsch-Gordan formula that
$\Lambda^k_\C M$ is a sum of representations of weight $\le k$.
By duality, the same is true about $\Lambda^{4n-k}_\C M$.
Denote by $\Lambda^k_+ M$ the maximal subrepresentation of weight
$k$ for $k \le 2n$ and of weight $4n-k$ for $k > 2n$ inside
$\Lambda^k_\C M$. Denote by $\Pi_+\colon \Lambda^* M\to \Lambda^*_+ M$
the equivariant projection onto the component of maximal weight.

Note that the Hodge decomposition with respect to an arbitrary
complex structure (we can pick $I$ without loss of generality)
is compatible with the weight decomposition:
$$
\Lambda^k_+ M = \bigoplus_{p+q=k} \Lambda^{p,q}_{I,+} M,
$$
where $\Lambda^{p,q}_{I,+} M = \Lambda^{p,q}_{I} M \cap \Lambda^k_+ M$.
This is actually a weight decomposition of an $SU(2)$-representation
for a special choice of Cartan subalgebra corresponding to $I$.
Therefore, $\Lambda^k_+ M$
is generated by $\Lambda^{k,0}_{I,+} M = \Lambda^{k,0}_I M$ as
an $SU(2)$-module. It follows that the summands $\Lambda^{p,q}_{I,+} M$
are of the same dimension and actually isomorphic
to $\Lambda^{p+q,0}_I M$. We will denote these isomorphisms by
$$
\mathcal{R}_{p,q}\colon \Lambda^{p+q,0}_I M \tilde{\to} \Lambda^{p,q}_{I,+} M,
$$
see \cite{_Verbitsky:qD_} and formula (\ref{_eqn_Rnn_}) below.

Given a hyperhermitian metric $g$ on $M$, we can consider the
corresponding K\"ahler forms $\omega_I$, $\omega_J$ and $\omega_K$.
It is easy to check that the form 
$$
\Omega_I = \omega_J + \sqrt{-1}\omega_K
$$
lies in $\Lambda^{2,0}_I M$. If $d\Omega_I = 0$, then the manifold
$M$ is called {\bf hyperk\"ahler}. This condition implies that
$g$ is K\"ahler with respect to any induced complex structure.
If $\Omega_I$ satisfies a strictly weaker condition $\partial\Omega_I = 0$,
then $M$ is called an {\bf HKT manifold}, and $g$ is called an
{\bf HKT metric} (HKT stands for HyperK\"ahler with Torsion;
see \cite{_Gra_Poon_} for a general introduction to HKT geometry).

We can reformulate the HKT condition as follows. Denote
by $d_+\colon \Lambda^k_+ M \to \Lambda^{k+1}_+ M$ the composition
$\Pi_+\circ d$ of the de Rham differential and the projection
onto the component of maximal weight. Note that
$\omega_I \in \Lambda^{1,1}_{I,+} M$. Then by Theorem 5.7 in
\cite{_Verbitsky:HKT_}, the condition $\partial\Omega_I = 0$
is equivalent to
\begin{equation}\label{_eqn_weight_}
d_+ \omega_I = 0.
\end{equation}
That is, for an HKT metric $g$ the exterior differential of its
K\"ahler form $d\omega_I$ has weight one. This observation will
be important for the proof of the main theorem.

\subsection[$SL(n, {\Bbb H})$-manifolds and holomorphic
Lagrangian calibration]{$SL(n, {\Bbb H})$-manifolds\\ and holomorphic
Lagrangian calibration}

As it was mentioned above, the holonomy of the Obata
connection on a hypercomplex manifold is a subgroup of
$GL(n, \H)$. Recall that this group can be defined as
follows. Consider a vector space $V$ of real dimension $4n$
with a quaternionic structure $I, J, K$. Then $GL(n, \H)$
consists of those linear transformations of $V$ that commute
with $I$, $J$ and $K$. Consider the Hodge decomposition
$V_\C = V\otimes \C = V^{1,0}_I\oplus V^{0,1}_I$.
The action of $GL(n, \H)$ preserves this decomposition,
hence it induces an action on all exterior powers of $V^{1,0}_I$.
Denote by $SL(n, \H)$ the subgroup of $GL(n, \H)$ consisting
of those elements that act identically on $V^{2n,0}_I$.

Consider a hypercomplex manifold $(M, I, J, K)$. If the
holonomy of the Obata connection is contained in $SL(n, \H)$,
then we call $M$ an {\bf $SL(n, \H)$-manifold}.
It follows almost immediately from the definition (see
\cite{_Verbitsky:canoni_}, Claim 1.1), that on an 
$SL(n, \H)$-manifold the canonical bundle $\Lambda^{2n,0}_{I} M$
is flat with respect to the Obata connection. It also
follows (\cite{_Verbitsky:canoni_}, Claim 1.2), that any
parallel section of $\Lambda^{2n,0}_{I} M$ is holomorphic.
Moreover, if the manifold $(M,I,J,K)$ is HKT and compact,
the condition $\Hol(M)\subset SL(n, {\Bbb H})$ is 
equivalent to holomorphic triviality of the canonical bundle
(\cite{_Verbitsky:canoni_}, Theorem 2.3).

Denote by $\Phi_I \in \Lambda^{2n,0}_{I} M$ a parallel section
which we call {\bf a holomorphic volume form}. Note that the
operator $J$ defines a real structure (that is, a 
complex-antilinear involution) on $\Lambda^{2n,0}_{I} M$ by
$\eta\mapsto \overline{J\eta}$. We can always assume
that $\Phi_I$ is real with respect to this structure.

Given an $SL(n, \H)$-manifold $M$ with a hyperhermitian
metric $g$, denote by $\Phi_I$ the holomorphic volume form
and by $\Omega_I = \omega_J+\sqrt{-1}\omega_K$ the 
$(2,0)$-form associated with the metric. 
In \cite{_Gra_Verb_}, a number of calibrations (in the
sense of \cite{_Harvey_Lawson:Calibrated_}) on 
$SL(n, \H)$-manifolds were constructed. In particular,
it was shown that the form
\begin{equation}\label{_eqn_Psi_}
\Psi = (-\sqrt{-1})^n\mathcal{R}_{n,n}(\Phi_I) \in \Lambda^{n,n}_{I,+} M
\end{equation}
is a calibration with respect to a properly rescaled metric.
It was shown that this form calibrates $\Omega_I$-Lagrangian
subvarieties of $M$. Recall, that a complex (with respect to
the complex structure $I$) subvariety $N\subset M$ of complex
dimension $n$ is called {\bf $\Omega_I$-Lagrangian} if the
restriction of $\Omega_I$ to the smooth part of $N$ vanishes.
This condition is equivalent to the following:
$T_x N$ is orthogonal to $J(T_x N)$ with respect to the metric
$g$ at any smooth point $x\in N$. To see this, observe that
for any vector fields $X, Y\in T^{1,0}_I M$ we have
$\Omega_I(X, Y) = 2 g(JX, Y)$.

We will not need the construction from \cite{_Gra_Verb_} in its
full generality, but only some basic properties of the form
$\Psi$. For the convenience of the reader we will formulate
these properties in the following lemma and give a short proof,
independent from \cite{_Gra_Verb_}.

\hfill

\lemma\label{_Psi_Lemma_}
Let $(M, I, J, K)$ be an $SL(n, \H)$-manifold of real dimension~$4n$,
$\Phi_I\in \Lambda^{2n,0}_I M$ a holomorphic volume form
and $\Psi = (-\sqrt{-1})^n\mathcal{R}_{n,n}(\Phi_I) \in \Lambda^{n,n}_{I,+} M$.
Let $N\subset M$ be an $I$-complex subvariety, $\dim_\C N = n$,
with $T_x N \cap J(T_x N) = 0$ at all smooth points $x\in N$.
Then:
\begin{enumerate}
\item $d\Psi = 0$,
\item $\Psi|_{N}$ is a strictly positive volume form on the smooth part of $N$.
\end{enumerate}

{\bf Proof:}
{1.} As it was mentioned above (see also formula (\ref{_eqn_Rnn_}) below),
the operator $\mathcal{R}_{n,n}$
is induced by the $SU(2)$-action on the exterior algebra of $M$. Since
the Obata connection $\nabla$ preserves the hypercomplex structure,
it commutes with the $SU(2)$-action and with the operator
$\mathcal{R}_{n,n}$. It was explained above that
$\nabla \Phi_I = 0$, thus $\nabla \Psi = 0$ also. The Obata
connection is torsion-free, hence $\nabla \Psi = 0$ implies $d\Psi = 0$.

{2.} Let $x\in N$ be a smooth point. By the assumptions of the
lemma, we can choose a basis $e_1,\ldots,e_n$ of $T^{1,0}_{I,x}N$,
such that $e_1,\ldots,e_n,J\overline{e}_1,\ldots,J\overline{e}_n$
will form a basis of $T^{1,0}_{I,x} M$. To show that
$\Psi|_{N}$ is strictly positive, we have to evaluate
the form $\Psi$ on the polyvector
$\xi = (-\sqrt{-1})^n e_1\wedge\overline{e}_1\wedge\cdots
\wedge e_n\wedge\overline{e}_n\in \Lambda^{n,n}_I (T_x M)$.

We will use the following explicit description of the operator $\mathcal{R}_{n,n}$.
Consider the operators: $\mathcal{H} = -\sqrt{-1}I$,
$\mathcal{X} = \frac{1}{2}
(\sqrt{-1}K - J)$ and $\mathcal{Y} = \frac 1 2 (\sqrt{-1}K + J)$
acting on $\Lambda^1_\C M$.
It is straightforward to check that these operators satisfy
the standard $\mathfrak{sl}_2(\C)$ commutator relations: $[\mathcal{X},\mathcal{Y}]=\mathcal{H}$,
$[\mathcal{H},\mathcal{X}]=2\mathcal{X}$, $[\mathcal{H},\mathcal{Y}]=-2\mathcal{Y}$.
Also observe that
$$
\mathcal{H}|_{\Lambda^{1,0}_{I} M} = 1,\hspace{10pt}
\mathcal{X}|_{\Lambda^{1,0}_{I} M} = 0,\hspace{10pt}
\mathcal{Y}|_{\Lambda^{1,0}_{I} M} = J,
$$
$$
\mathcal{H}|_{\Lambda^{0,1}_{I} M} = -1,\hspace{10pt}
\mathcal{X}|_{\Lambda^{0,1}_{I} M} = -J,\hspace{10pt}
\mathcal{Y}|_{\Lambda^{0,1}_{I} M} = 0.
$$

We can extend $\mathcal{H}$, $\mathcal{X}$, $\mathcal{Y}$ as derivations
to the whole exterior algebra. Then we have 
\begin{equation}\label{_eqn_Rnn_}
\mathcal{R}_{n,n} = \mathcal{Y}^n,
\end{equation}
and
$$
\langle \Psi, \xi\rangle = (-\sqrt{-1})^n\langle \mathcal{Y}^n\Phi_I, \xi \rangle
= (\sqrt{-1})^n\langle \Phi_I, \mathcal{Y}^n\xi \rangle.
$$
Observe, that $\mathcal{Y}|_{T^{1,0}_{I} M} = 0$, $\mathcal{Y}|_{T^{0,1}_{I} M} = J$,
and since $\mathcal{Y}$ acts as a derivation, we have
$$
\mathcal{Y}^n\xi = n!(-\sqrt{-1})^n e_1\wedge J\overline{e}_1\wedge\cdots \wedge e_n\wedge J\overline{e}_n.
$$
The volume form is given by $\Phi_I = a\, e^*_1\wedge J\overline{e^*_1}\wedge\cdots \wedge e^*_n\wedge J\overline{e^*_n}$
for some positive real number $a$, and we see that $\langle\Psi, \xi\rangle > 0$.
\endproof


\section[Holomorphic Lagrangian fibrations on $SL(n, {\Bbb H})$-manifolds]
{Holomorphic Lagrangian fibrations\\  on $SL(n, {\Bbb H})$-manifolds}


\subsection{HKT structures and Lagrangian fibrations}

\definition
\label{_holo_lagra_Definition_}
A {\bf holomorphic Lagrangian subvariety} 
of an $SL(n, {\Bbb H})$-manifold
is a subvariety calibrated by the form $\Psi$ of 
\ref{_Psi_Lemma_}.

\hfill

\definition
Let $M$ be an $SL(n, {\Bbb H})$-manifold, and $\phi\colon (M,I) \to X$
a smooth holomorphic fibration. It is called {\bf a holomorphic
Lagrangian fibration} if all its fibers are holomorphic Lagrangian
subvarieties. We will say that the fibration $\phi\colon (M,I) \to X$ is
{\bf smooth} if $\phi$ is a submersion (hence all the fibers are smooth).

\hfill

\remark
Suppose that we are given an arbitrary holomorphic fibration 
$\phi\colon (M,I) \to X$ with $M$ compact and $\dim_\C X = {1\over 2}\dim_\C M$.
Denote by $F_x$ the fiber of $\phi$ over $x\in X$. Consider the
set $U=\{x\in X\colon J(TF_x)\cap TF_x = 0\}$. It is clear that
$U$ is open and all the fibers over $U$ can be made Lagrangian with
respect to some hyperhermitian metric. Namely, the condition
$J(TF_x)\cap TF_x = 0$ implies that each fiber of the vector bundle
$J(TF_x)$ projects onto $T_x X$ by $\phi$, so that we can
lift any Hermitian metric from $U$ to $J(TF_x)$. Then we can uniquely 
extend it to a hyperhermitian metric on $\phi^{-1}(U)$ with $J(TF_x)$
orthogonal to $TF_x$. This means that the set of Lagrangian fibers
of $\phi$ is open in $M$. 

\hfill

The main result of this paper is the following theorem.

\hfill

\theorem\label{_HKT=>Kahler_Theorem_}
Let $M$ be a compact $SL(n, {\Bbb H})$-manifold, and $\phi\colon (M,I) \to X$
a smooth holomorphic Lagrangian fibration. Assume that $M$ admits an
HKT-structure. Then $X$ is K\"ahler.

\hfill

{\bf Proof:} Let $\omega_I\in \Lambda^{1,1}_I M$
be a K\"ahler form, and $\Psi$ the holomorphic Lagrangian calibration
defined by (\ref{_eqn_Psi_}). Consider the form 
$$
\Theta = \Psi\wedge \omega_I \in \Lambda^{n+1,n+1}_I M.
$$
By \ref{_Psi_Lemma_}, $\Psi$ is closed, so
we have $d\Theta = \Psi\wedge d\omega_I \in \Lambda^{2n+3} M$.
It follows from (\ref{_eqn_weight_}) that $d\omega_I$ has weight one,
hence Clebsch-Gordan formula implies that 
$d\Theta$ is a sum of weight $2n+1$ and weight $2n-1$ components.
But the weight decomposition of $\Lambda^{2n+3} M$ has components
only up to weight $2n-3$, so $d\Theta$ has to vanish.

We will obtain a K\"ahler metric on $X$ as a push-forward
$\pi_*\Theta$, which is a $(1,1)$-current on $X$.
We have to show that this current is strictly positive
and smooth.

Consider any Hermitian metric on $X$ and denote by $\eta$
the corresponding $(1,1)$-form. It is always possible to
rescale $\eta$ so that it satisfies the condition $\pi^*\eta \le \omega_I$.
In this case we have $\Theta = \Psi\wedge\omega_I \ge \Psi\wedge\pi^*\eta$.
The inequality is preserved under taking push-forwards, so
we have $\pi_*\Theta \ge \pi_*(\Psi\wedge\pi^*\eta)$.

For any test-form $\alpha\in \Lambda^{n-1,n-1}_I X$ we have by definition
\begin{eqnarray}
\int_X \pi_*(\Psi\wedge\pi^*\eta)\wedge \alpha = \int_M \Psi\wedge\pi^*\eta\wedge \pi^*\alpha\nonumber\\
= \int_M \Psi\wedge\pi^*(\eta \wedge \alpha) = \int_X \pi_*\Psi\wedge \eta \wedge \alpha.\nonumber
\end{eqnarray}

Thus, we have the inequality $\pi_*\Theta \ge \pi_*(\Psi)\wedge \eta$ where
$\pi_*(\Psi)$ is a closed \hbox{$0$-current} on $X$, that is a constant.
To see that this constant is positive, take $\alpha = \eta^{n-1}$
and observe that $\Psi\wedge \pi^*(\eta^n)$ is a smooth positive form
of top degree on $M$. It suffices to check that this form is non-zero at
some point of $M$. Pick any non-critical point $x$ of $\pi$ and decompose the
tangent space $T_x M = V_0 \oplus V_1$, where $V_0$ is tangent to the fiber
and $V_1$ is any complementary subspace. Since $\pi^*\eta^n|_{V_0} = 0$
and $\pi^*\eta^n|_{V_1}$ is a volume form on $V_1$, we should
check that $\Psi|_{V_0}$ is non-zero. But this follows from the second
statement in \ref{_Psi_Lemma_}, since the fiber is Lagrangian.

It remains to check that $\pi_*\Theta$ is smooth. This follows from
smoothness of the fibration $\phi\colon (M,I) \to X$. Since $\phi$
is a submersion, we can choose a splitting $TM = \mathcal{V}\oplus \mathcal{H}$,
where $\mathcal{V}$ is a subbundle, tangent to the fibers of $\phi$.
Using this splitting we can lift vector field from $X$ to $M$.
Denote by $F_x$ the fiber $\phi^{-1}(x)$ over a point $x\in X$.
Consider a $(1,1)$-form $\theta$ on $X$, which is given by
$$
\theta(\xi,\overline{\eta})_x = \int_{F_x} (\pi^*\xi)\lrcorner\, (\pi^*\overline\eta)\lrcorner\, \Theta,
$$
where $\xi, \eta \in T^{1,0} X$ and $\pi^*\xi$, $\pi^*\eta$ denote
the the lifts of these vector fields to $\mathcal{H}$.
Note, that $\theta$ does not actually depend on the choice of $\mathcal{H}$,
since convolution of $\Theta$ with a vertical vector field from $\mathcal{V}$
gives a differential form that restricts trivially to the fiber.
The form $\theta$ coincides with $\pi_*\Theta$ as a current, thanks to Fubini
theorem. But it is clear from the definition that $\theta$ is smooth.
This completes the proof.
\endproof

\hfill

\remark
Note that for a non-smooth holomorphic Lagrangian fibration
$\phi\colon (M,I) \to X$ the above proof shows that 
$\pi_*(\Theta \wedge \omega_I)$ is a K\"ahler current. Therefore,
$X$ is in Fujiki class $\mathcal{C}$, whenever $M$ admits an HKT-structure
(see \cite{_Dem_Paun_}).

\subsection{Examples of Lagrangian fibrations}
\label{_exa_Lagra_Subsection_}

In this subsection we construct a class of hypercomplex
$SL(n, \H)$-manifolds that do not admit an HKT metric.
The idea, which was also used in \cite{_Konts_Soib_},
is to consider a torus fibration over an affine base.

Let $(X, I)$ be a complex manifold, $\dim_\C X = n$.
We call $X$ {\bf an affine complex manifold} if
there exists a torsion-free connection
$D\colon TX\to TX\otimes \Lambda^1 X$ which is flat
and preserves the complex structure: $DI = 0$.

Fix a point $x\in X$ and consider the holonomy
group $\mathrm{Hol}_x(D)\subset GL(T_x X)$. We call
$X$ an affine manifold {\bf with integer monodromy}
if there exists a lattice $\Lambda_x \subset T_x X$
which is preserved by holonomy: 
$\mathrm{Hol}_x(D)\Lambda_x = \Lambda_x$.
In this case, we can construct a subset $\Lambda\subset TX$
in the total space of the tangent bundle, obtained as
parallel translation of $\Lambda_x$.
For any point $y\in X$ the intersection $\Lambda\cap T_y X$
is a lattice in $T_y X$.

Let $M = TX / \Lambda$, that is, a manifold obtained as
a fiberwise quotient of $TX$. We will introduce a pair of
anticommuting complex structures on $TX$ that will descend to $M$.
The connection $D$ defines a splitting 
\begin{equation}\label{_eqn_split_}
T(TX) = \mathcal{V}\oplus \mathcal{H}
\end{equation}
into a direct sum of the vertical subbundle $\mathcal{V}$
which is tangent to the fibers of the projection $TX\to X$,
and a horizontal complement $\mathcal{H}$.
Note that for any point $(x,v)\in TX$ with $v\in T_x X$ we
have natural isomorphisms $\mathcal{H}_{(x,v)}\simeq T_x X \simeq \mathcal{V}_{(x,v)}$,
thus we can identify the components of the splitting (\ref{_eqn_split_}).
This also shows that the complex structure $I$ acts
naturally on both $\mathcal{H}$ and $\mathcal{V}$.
With these observations in mind, define a pair of
operators from $\End(T(TM))$:
$$
\mathcal{I} = 
\begin{pmatrix}
-I & 0\\ 0 & I
\end{pmatrix},\quad
\mathcal{J} =
\begin{pmatrix}
0 & \Id\\ -\Id & 0
\end{pmatrix},
$$
where the block-matrix form corresponds to the decomposition (\ref{_eqn_split_}).
It is clear that these operators define a pair of anticommuting
almost-complex structures. We need to check that these structures
are integrable. This can be done locally.

Using the fact that the connection $D$ is flat and $I$ is parallel,
we can choose a parallel local frame in $TX$ of the form 
$e_1, Ie_1, \ldots e_n, Ie_n$. Since $D$ is torsion-free,
these vector fields commute, hence they define a local coordinate
system on $X$. Since the sections $e_i$ and $Ie_i$ are flat,
they are tangent to the subbundle $\mathcal{H}$. Therefore, they also
define a local coordinate system in the total space of $TX$
in which the complex structures $\mathcal{I}$ and $\mathcal{J}$
act as standard complex structures of $\H^n$. Thus, they are
integrable.

This construction gives a hypercomplex structure on the total
space of $TX$. Observe that this structure is invariant with respect
to translations along the fiber. This implies that the hypercomplex structure
descends to $M$. We call the manifold $M$ constructed this way
{\bf quaternionic double} of an affine complex manifold $X$.

The hypercomplex structure on the quaternionic double $M$ is locally
modeled on $\H^n$ by construction. Thus, the Obata connection
on this manifold is flat (and it is actually induced from the
flat connection $D$ on $X$). Since the connection $D$ preserves
the lattice in the tangent bundle, it also preserves the
holomorphic volume form. To see this, we can again choose a local
frame $e_1, Ie_1, \ldots e_n, Ie_n$ as above. We can assume that
the elements of this frame form a local basis for the lattice $\Lambda$ preserved
by $D$. Then the volume form
$(e_1^* - \sqrt{-1}I e_1^*)\wedge\cdots \wedge(e_n^* - \sqrt{-1}I e_n^*)$
is well-defined globally and parallel with respect to $D$,
hence holomorphic. It follows that the Obata connection
on $M$ also preserves an $I$-holomorphic volume form.
We conclude that the manifold $M$ is an $SL(n, \H)$-manifold.

\hfill

\theorem
Let $M$ be a quaternionic double of an affine complex manifold $X$.
If $X$ is not K\"ahler, then $M$ does not admit an HKT metric. 

{\bf Proof:} We can choose a Hermitian metric $g$ on $X$.
Then we can define a hyperhermitian metric $h=g\oplus g$
on $M$, using the decomposition (\ref{_eqn_split_}).
With respect to this metric the fibers of the projection
$M\to X$ are Lagrangian, hence we can apply \ref{_HKT=>Kahler_Theorem_}. 
\endproof

\hfill

We conclude this section by some examples of affine complex
manifolds with integer monodromy. We remark that the
manifolds in the examples will be non-K\"ahler (actually
any complex affine manifold admitting a K\"ahler metric
is a quotient of a torus, as follows from the 
Calabi-Yau theorem and the Bieberbach's solution of Hilbert's 18-th
problem; see e.g. \cite{_Besse:Einst_Manifo_}). 
Therefore, the corresponding quaternionic 
doubles possess no HKT metrics.

\hfill

\example
Let $N$ be a three-dimensional real Heisenberg group,
that is the group of upper-triangular unipotent $3\times 3$
matrices with real elements.
Consider the nilpotent Lie group $G = N\times \R$. Then
$G$ is diffeomorphic to $\C^2$ and one can check
(see \cite{_Hasegawa_}, Example 2) that the structure
of the Lie group on $\C^2$ can be explicitly given by
$$
(w_1,w_2)\cdot(z_1,z_2) = (w_1+z_1,w_2-\sqrt{-1}\overline{w_1}z_1+z_2).
$$ 
It is clear from this formula, that left translations
are holomorphic affine transformations of $C^2$, so $G$
has a left-invariant complex structure
which is preserved by a flat torsion-free connection.
We can choose a lattice in $\C^2$, say $\Lambda = \Z[\sqrt{-1}]^2$.
Then the manifold $X = \Lambda\backslash G$
(which is called {\bf a primary Kodaira surface}) will be an
affine complex manifold with integer monodromy.

\hfill

\example
Let $G$ be a three-dimensional complex Heisenberg group.
It can be described (see \cite{_Hasegawa_}) as $\C^3$
with multiplication given by
$$
(w_1,w_2,w_3)\cdot(z_1,z_2,z_3) = (w_1+z_1,w_2+z_2,w_3+z_3+w_1 z_2).
$$
We can choose a lattice in $G$, for example $\Lambda=\Z[\sqrt{-1}]^3$
and take a quotient $X = \Lambda\backslash G$. This is
an affine complex manifold with integer monodromy (which is called 
{\bf an Iwasawa manifold}).

\hfill

\hfill

\hfill

\small{

\noindent

{\sc Misha Verbitsky}\\
{\sc Laboratory of Algebraic Geometry, \\
National Research University HSE,\\
7 Vavilova Str. Moscow, Russia, 117312}\\
{\tt verbit@maths.gla.ac.uk, \ \  verbit@mccme.ru\\}

{\sc Andrey Soldatenkov\\
{\sc Laboratory of Algebraic Geometry,\\
National Research University HSE,\\
7 Vavilova Str. Moscow, Russia, 117312}\\
{\tt aosoldatenkov@gmail.com}
 }}


{\small
\begin{thebibliography}{AV1}







\bibitem[BD]{_Barberis+Dotti:1996_}
M. L. Barberis, I. Dotti, 
{\em Hypercomplex structures on a class of solvable Lie groups,} 
Quarterly Journal of Mathematics Oxford (2),   47 (1996), 389-404. 

 

\bibitem[BDV]{_BDV:nilmanifolds_}
Maria L. Barberis, Isabel G. Dotti, Misha Verbitsky, {\em
  Canonical bundles of complex nilmanifolds, with
applications to hypercomplex geometry}, arXiv:0712.3863, 22 pages.


\bibitem[Bes]{_Besse:Einst_Manifo_}
Besse, A., {\em Einstein Manifolds}, Springer-Verlag, New York (1987)


\bibitem[Bo]{_Boyer_} 
Boyer, Charles P.
{\em A note on hyper-Hermitian four-manifolds}.
Proc. Amer. Math. Soc. 102 (1988), no. 1, 157--164. 



\bibitem[DP]{_Dem_Paun_}
Demailly, J.-P., P\u{a}un, M., {\em Numerical characterization of the K\"ahler cone of a compact
K\"ahler manifold}, Annals of Math, 159 (2004), 1247--1274.

\bibitem[FG]{_Fino_Gra_}
Fino, A.,  Grantcharov, G., {\em On some properties of the
  manifolds with skew-symmetric torsion and holonomy
SU(n) and Sp(n)}, math.DG/0302358, Adv. Math. 189 (2004), no. 2, 439--450.


\bibitem[GP]{_Gra_Poon_}
Grantcharov, G., Poon, Y. S., {\em Geometry of
  hyper-K\"ahler connections with torsion},
 math.DG/9908015,
Comm. Math. Phys. 213 (2000), no. 1, 19--37.


\bibitem[GV]{_Gra_Verb_}
Grantcharov, G., Verbitsky, M., {\em Calibrations in hyperk\"ahler geometry},
arXiv:1009.1178, 32 pages.


\bibitem[GP1]{_GP1_} 
Gutowski, J., Papadopoulos, G. 
{\em The Dynamics of Very Special Black Holes},
Phys. Lett. B472 (2000) 45-53


\bibitem[GP2]{_GP2_} 
Gutowski, J., Papadopoulos, G. {\em
The Moduli Spaces of Worldvolume Brane Solitons},
 Phys. Lett. B432 (1998) 97-102.


\bibitem[GPS]{_GPS_} 
G.W. Gibbons, G. Papadopoulos, K.S. Stelle,
{\em HKT and OKT Geometries on Soliton Black Hole Moduli Spaces},
Nucl.Phys. B508 (1997) 623-658



\bibitem[HL]{_Harvey_Lawson:Calibrated_}
Harvey, R., Lawson, B., {\em Calibrated geometries},  Acta
Math.  148 (1982), 47-157.


\bibitem[Ha]{_Hasegawa_}
Hasegawa, K., {\em Complex and K\"ahler structures on compact homogeneous manifolds --
their existence, classification and moduli problem},
arXiv:0804.0738, 14 pages.




\bibitem[HP]{_Howe_Papado_}
P.S. Howe, G. Papadopoulos,  
{\em Twistor spaces for hyper-K\"ahler manifolds with torsion} Phys. Lett. B 379
(1996), no. 1-4, 80--86.




\bibitem[J]{_Joyce_}
D. Joyce,  {\em Compact hypercomplex and quaternionic
  manifolds}, J. Differential Geom. {\bf 35} (1992) no. 3,
743-761


\bibitem[KRS]{_KRS_} A.Kapustin, 
L. Rosansky, N. Saulina {\em Three dimensional topological 
field theory and symplectic algebraic geometry I}, 
Nucl Phys. B {\bf 816} (2009), no. 3, 295-355.


\bibitem[KS]{_Konts_Soib_}
Kontsevich, M., Soibelman, Y., {\em Homological mirror symmetry and torus fibrations},
arXiv:math/0011041, 66 pages. 






\bibitem[Ob]{_Obata_}
Obata, M., {\em Affine connections on manifolds with
  almost complex, quaternionic or Hermitian structure},
Jap. J. Math., 26 (1955), 43-79.


\bibitem[SSTV] {_SSTvP_}
Ph. Spindel, A. Sevrin, W. Troost, A. Van Proeyen 
{\em Extended
supersymmetric $\sigma$-models on group manifolds}, Nucl. Phys. B308
(1988) 662-698.


\bibitem[Saw]{_Sawon_}
Sawon, J.
{\em Abelian fibred holomorphic symplectic
  manifolds}, Turkish
Jour. Math. 27 (2003), no. 1, 197-230, math.AG/0404362.


\bibitem[Sol]{_Soldatenkov:SU(3)_} 
Andrey Soldatenkov,
{\em Holonomy of the Obata connection on $SU(3)$}, 
arXiv:1104.2085, Int. Math. Res. Notices 
(2012), Vol. 2012 (15), 3483-3497.

\bibitem[SV]{_SV:trianalytic_} 
Andrey Soldatenkov, Misha Verbitsky
{\em Subvarieties of hypercomplex manifolds with holonomy
  in $SL(n,{\Bbb H})$}
Journal of Geometry and Physics, Volume 62,
Issue 11, November 2012, Pages 2234-2240 



\bibitem[V0]{_Verbitsky:trianalyt_}
Verbitsky M., {\em Tri-analytic subvarieties of hyper-Kaehler manifolds,} 
also known as {\em Hyperk\"ahler embeddings and holomorphic 
symplectic geometry II}, GAFA {\bf 5} no. 1 (1995), 92-104, 
alg-geom/9403006.


\bibitem[V1]{_Verbitsky:HKT_}
Verbitsky, M., {\em Hyperk\"ahler manifolds with torsion,
  supersymmetry and Hodge theory}, math.AG/0112215,
Asian J. Math. Vol. 6, No. 4, pp. 679-712 (2002).


\bibitem[V2]{_Verbitsky:qD_}
Verbitsky, M., {\em Quaternionic Dolbeault complex and
  vanishing theorems on hyperkahler manifolds},
 Compos. Math. 143 (2007), no. 6, 1576--1592, math/0604303.


\bibitem[V3]{_Verbitsky:canoni_}
Verbitsky, M., {\em Hypercomplex manifolds with trivial canonical bundle
 and their holonomy}, arXiv:math/0406537,
``Moscow Seminar on Mathematical Physics, II'', American Mathematical Society Translations, {\bf 2}, 221 (2007).

\bibitem [V4]{_V:hypercomplex_Kahler_}
Verbitsky, M., {\em 
Hypercomplex structures on K\"ahler manifolds},
math.AG:0406390, 10 pages, GAFA 15 (2005), no. 6,
pp. 1275-1283.






\end{thebibliography}
}
\end{document}